\documentclass[a4paper,12pt]{article}
\usepackage{latexsym}
\usepackage{amsmath}
\usepackage{amssymb}
\usepackage{amscd}
\usepackage{array}
\usepackage{comment}
\usepackage[all]{xy}
\usepackage{amsthm}
\usepackage{amsfonts}
\usepackage{enumerate}
\newtheorem{theorem}{Theorem}[]
\newtheorem{proposition}[theorem]{Proposition}

\theoremstyle{definition}

\newcommand{\HH}{\mathcal H}

\newcommand{\Gal}{\mathrm{Gal}}

\newcommand{\Hol}{\mathrm{Hol}}

\newcommand{\Sym}{\operatorname{Sym}}

\newcommand{\SL}{\mathrm{SL}}

\newcommand{\End}{\operatorname{End}}

\newcommand{\Aut}{\operatorname{Aut}}

\newcommand{\wK} {{\widetilde{K}}}

\title{Hopf Galois structures on symmetric and alternating extensions}
\author{Teresa Crespo, Anna Rio and Montserrat Vela}
\date{\today}

\begin{document}

\maketitle

\begin{abstract} By using our previous results on induced Hopf Galois structures and a recent result by Koch, Kohl, Truman and Underwood on normality, we determine which types of Hopf Galois structures occur on Galois extensions with Galois group isomorphic to alternating or symmetric groups. \let\thefootnote\relax\footnotetext{T. Crespo acknowledges support by grants MTM2015-66716-P (MINECO/FEDER, UE) and 2014SGR 206; A. Rio and M. Vela acknowledge support by grants MTM2015-66180R (MINECO/FEDER, UE) and 2014SGR 550.}



\end{abstract}

\section{Introduction}

A Hopf Galois structure on a finite extension of fields $K/k$ is a pair $(\HH,\mu)$, where $\HH$ is
a finite cocommutative $k$-Hopf algebra  and $\mu$ is a
Hopf action of $\HH$ on $K$, i.e a $k$-linear map $\mu: \HH \to
\End_k(K)$ giving $K$ a left $\HH$-module algebra structure and inducing a bijection $K\otimes_k \HH\to\End_k(K)$.
Hopf Galois structures were introduced by Chase and Sweedler in \cite{C-S}.

In Hopf Galois theory one has the following Galois correspondence theorem.

\begin{theorem}[\cite{C-S} Theorem 7.6] Let $(\HH,\mu)$ be a Hopf Galois structure on the field extension $L/K$.
For a $K$-sub-Hopf algebra $\HH'$ of $\HH$ we define
$$
L^{\HH'}=\{x\in L \mid \mu(h)(x)=\varepsilon(h)\cdot x \mbox{ for all } h\in \HH'\},
$$
where $\varepsilon$ is the counity of $\HH$.
Then, $L^{\HH'}$ is a subfield of $L$, containing $K$, and
$$
\begin{array}{rcl}
{\mathcal F}_{\HH}:\{\HH'\subseteq \HH \mbox{ sub-Hopf algebra}\}&\longrightarrow&\{\mbox{Fields }E\mid K\subseteq E\subseteq L\}\\
\HH'&\to &L^{\HH'}
\end{array}
$$
is injective and inclusion reversing.
\end{theorem}

For separable field extensions, Greither and
Pareigis \cite{G-P} give the following group-theoretic
equivalent condition to the existence of a Hopf Galois structure.

\begin{theorem}
Let $K/k$ be a separable field extension of degree $n$, $\wK$ its Galois closure, $G=\Gal(\wK/k), G'=\Gal(\wK/K)$. Then there is a bijective correspondence
between the set of Hopf Galois structures on $K/k$ and the set of
regular subgroups $N$ of the symmetric group $S_n$ normalized by $\lambda (G)$, where
$\lambda:G \rightarrow S_n$ is the morphism given by the action of
$G$ on the left cosets $G/G'$.
\end{theorem}

For a given Hopf Galois structure on $K/k$, we will refer to the isomorphism class of the corresponding group $N$ as the type of the Hopf Galois
structure. The Hopf algebra $\HH$ corresponding to a regular subgroup $N$ of $S_n$ normalized by $\lambda (G)$ is the Hopf subalgebra $\wK[N]^G$ of the group algebra $\wK[N]$ fixed under the action of $G$, where $G$ acts on $\wK$ by $k$-automorphisms and on $N$ by conjugation through $\lambda$. It is known that the Hopf subalgebras of $\wK[N]^G$ are in 1-to-1 correspondence with the subgroups of $N$ stable under the action of $G$ (see e.g. \cite{CRV} Proposition 2.2). For $N'$ a $G$-stable subgroup of $N$, we will denote by $K^{N'}$ the subfield $K^{\HH'}$ of $K$ fixed by the Hopf subalgebra $\HH'$ of $\HH$ corresponding to $N'$ and refer to it as fixed by $N'$.

Childs \cite{Ch1} gives an equivalent  condition to the existence of a Hopf Galois structure introducing the holomorph of the regular subgroup $N$ of $S_n$. We state the more precise formulation of this result due to Byott \cite{B} (see also \cite{Ch2} Theorem 7.3).

\begin{theorem}\label{theoB} Let $G$ be a finite group, $G'\subset G$ a subgroup and $\lambda:G\to \Sym(G/G')$ the morphism given by the action of
$G$ on the left cosets $G/G'$.
Let $N$ be a group of
order $[G:G']$ with identity element $e_N$. Then there is a
bijection between
$$
{\cal N}=\{\alpha:N\hookrightarrow \Sym(G/G') \mbox{ such that
}\alpha (N)\mbox{ is regular}\}
$$
and
$$
{\cal G}=\{\beta:G\hookrightarrow \Sym(N) \mbox{ such that }\beta
(G')\mbox{ is the stabilizer of } e_N\}.
$$
Under this bijection, if $\alpha\in {\cal N}$ corresponds to
$\beta\in {\cal G}$, then $\alpha(N)$ is normalized by
$\lambda(G)$ if and only if $\beta(G)$ is contained in the
holomorph $\Hol(N)$ of $N$.
\end{theorem}

\section{Main result}
We will apply Theorem 2.9 in \cite{KKTU} to discard some types of Hopf Galois structures on Galois extensions with given Galois group. The setting will be the following. Let $G$ be a group of order $n$ and $G'$ a subgroup of $G$ of index $d$ such that no nontrivial subgroup of $G'$ is normal in $G$. Let $N$ be a group of order equal to $n$ having a unique conjugation class of subgroups of index $d$ with length 1. With these hypothesis, if $K/k$ is a Galois extension  with Galois group $G$ and $F:=K^{G'}$, then $F:=K^{N'}$, for $N'$ the normal subgroup of $N$ of index $d$. If we know that a separable extension of degree $d$ having normal closure with Galois group $G$ has no Hopf Galois structure of type $N/N'$, we may conclude that a Galois extension with Galois group $G$ has no Hopf Galois structure of type $N$.
We will use Theorem 3 in \cite{CRV3} to prove that a certain type of Hopf Galois structure does occur on a Galois extension with given Galois group. Let $K/k$ be  a Galois extension  with Galois group $G=H\rtimes G'$ and let $F:=K^{G'}$. Then if $F/k$ has a Hopf Galois structure of type $N_1$ and $K/F$ has a Hopf Galois structure of type $N_2$, the extension $K/k$ has a Hopf Galois structure of type $N_1\times N_2$.

\subsection{Galois extensions with Galois group $A_4$ or $S_4$}

Let us denote by $D_{2n}$ the dihedral group of order $2n$ and by $Dic_n$ the dicyclic group of order $4n$, that is

$$D_{2n}=\langle r,s | r^n=1, s^2=1, srs=r^{-1}\rangle, \quad Dic_n=\langle a,x | a^{2n}=1, x^2=a^n, xax^{-1}=a^{-1}\rangle.$$

Let us assume that $K/k$ is Galois with group $G=A_4$, the alternating group.
We analize the five possible types of Hopf Galois structures: the alternating group
$A_4$, the dicyclic group $Dic_3=C_3\rtimes C_4$, the cyclic group $C_{12}=C_3\times C_4$,  the dihedral group $D_{12}=C_3\rtimes V_4$ and the direct product
$C_3\times V_4$.

The classical Galois structure gives a Hopf Galois structure of type $A_4$. On the other hand, since
$A_4=V_4\rtimes C_3$, a quartic extension having Galois closure $A_4$ is Hopf Galois of type $V_4$, hence,
by \cite{CRV3}, theorem 3, we get induced Hopf Galois structures of type $C_3\times V_4$. Finally, since
$\Hol(Dic_3)=\Hol(D_{12})$ (see \cite{K}, Proposition 2.1) either both types of Hopf Galois structures arise or none of them does.
We are left with cyclic and dicyclic types.

In both cases,  we consider $N'$ the cyclic subgroup of order 3, the 3-Sylow subgroup, and we have
$N/N'\simeq C_4$. Then, the corresponding fixed field $F$ gives a quartic extension with Galois closure $K$.
Since $\Hol(C_4)$ has order 8, it can not contain $G$,
and this extension $F/k$ can not have Hopf Galois structures of type $C_4$.
This proves that $K/k$ has neither cyclic nor dicyclic (or dihedral) Hopf Galois structures.

Now let us assume that $K/k$ is Galois with group $G=S_4$, the symmetric group.
There are 15 isomorphism classes of groups of order 24. Hence there are 15 possible types for Hopf Galois structures on $K/k$.
For a group of order 24, the number $n_3$ of 3-Sylow subgroups may be 1 or 4. If $n_3=1$, the group is a semi-direct product $C_3\rtimes S$, where $S$ is a group of order 8, i.e. $S=C_8, C_4\times C_2, E_8=C_2\times C_2 \times C_2$, the dihedral group $D_8$ or the quaternion group $Q_8$. There are 12 groups of order 24 with $n_3=1$. These are precisely $C_3\rtimes C_8, C_{24}=C_3\times C_8, S_3\times C_4=C_3\rtimes (C_2\times C_4), Dic_{3}\times C_2=C_3\rtimes (C_2\times C_4), S_3\times V_4=C_3\rtimes (C_2\times C_2\times C_2), C_3\times (C_2\times C_2\times C_2), D_{24}=C_3\rtimes D_8,  C_3\rtimes_{\varphi} D_8$, where $\varphi:D_8 \rightarrow \Aut C_3$ has kernel $C_2\times C_2$, $C_3\times D_8, Dic_{6}= C_3\rtimes Q_8, C_3\times Q_8$.  We have $n_3=4$ for $S_4, SL(2,3)$ and $A_4\times C_2$.

Let us consider an intermediate field $F$ for the extension $K/k$ such that $[F:k]=8$. Then $F/k$ has Galois closure $K$ and by \cite{C-S}, $F/k$ has only Hopf Galois structures of type $E_8=C_2\times C_2\times C_2$. By \cite{KKTU} theorem 2.9, $K/k$ has no Galois structure of type $N$ if $N$ has a unique subgroup $N'$ of order 3 (then normal and $G$-stable) such that $N/N'$ is not isomorphic to $E_8$. This is the case for $N=C_3\rtimes C_8, C_{24}, S_3\times C_4, Dic_{3}\times C_2, D_{24},  C_3\rtimes_{\varphi} D_8, C_3\times D_8, Dic_{6}, C_3\times Q_8$.

Let us consider now the subfield $F$ of $K$ fixed by a transposition of $S_4$. Since $A_4$ is a normal complement of $\Gal(K/F)$ in $S_4$, the extension $F/k$ has a Hopf Galois structure of type $A_4$, hence by \cite{CRV3}, theorem 3,  $K/k$ has an induced Hopf Galois extension of type $A_4\times C_2$. Let us now take $F$ to be the subfield of $K$ fixed by a subgroup of $S_4$ isomorphic to $S_3$. Then $K/F$ is Galois with group $S_3$ and has a Hopf Galois structure of type $C_6$ (see \cite{C-S}). Since $F/k$ has a Hopf Galois structure of type $C_2\times C_2$, we obtain, again by \cite{CRV3}, theorem 3 and taking into account $S_4=V_4\rtimes S_3$, than $K/k$ has induced Hopf Galois structures of types $S_3\times V_4$ and $C_6\times V_4$. Finally, we check, using Magma, that $\Hol(\SL(2,3))$ has no subgroup isomorphic to $S_4$, hence $K/k$ has no Galois structure of type $\SL(2,3)$.

We have obtained the following result.

\begin{proposition}
Let $K/k$ be a Galois extension with Galois group $A_4$. Then, the only types of Hopf Galois structures on $K/k$ are $A_4$ and $C_3\times V_4$.
The classical Galois structure realizes type $A_4$ and a Hopf Galois structure of type $C_3\times V_4$ is induced by the classical Galois structure on $K/F$ and the Hopf Galois structure of type $V_4$ on $F/k$ for $F$ an intermediate field with $[K:F]=3$.

Let $K/k$ be a Galois extension with Galois group $S_4$. Then, the only types of Hopf Galois structures on $K/k$ are $S_4$
and the split ones $A_4\times C_2, S_3\times V_4$ and $C_6\times V_4$. The classical Galois structure realizes the first type and the remaining three are realized as induced structures.
\end{proposition}

\subsection{Galois extensions with Galois group $A_5$ or $S_5$}\label{5}

Let us assume that $K/k$ is Galois with Galois group $G=A_5$, the alternating group.
There are 13 possible types of Hopf Galois structures. If we take $N\ne A_5$ a group of order 60, then
$N$ has a unique $5-$Sylow subgroup that we can take for $N'$. Since none of the groups of order 12 has
holomorph of order divisible by 60, we know that $K^{N'}/k$ is not a Hopf Galois extension and therefore theorem 2.9 in \cite{KKTU} implies
than $N$ is not a Hopf Galois type for $K/k$.

Let us assume that $K/k$ is Galois with group $G=S_5$, the symmetric group.
Now there are 47 possible types of Hopf Galois structures.

The classical Galois structure gives a Hopf Galois structure of type $S_5$. On the other hand, since
$S_5=A_5\rtimes C_2$, again by \cite{CRV3}, theorem 3, we get induced Hopf Galois structures of type $A_5\times C_2$.
(Recall that an extension $F/k$ of degree 60 with Galois closure $K/k$ has an almost classical Hopf Galois structure of
type $A_5$, since $A_5$ is a normal complement of $\Gal(K/F)$ in $\Gal(K/k)\simeq S_5$.)

Checking on the remaining 45 types, we see that all but one, namely $N=\SL(2,5)$, have a normal $p-$Sylow subgroup.
For a given $N$ in that set, we choose $N'$ a normal $p-$Sylow subgroup. Therefore, $N'$ is a normal $G$-stable subgroup
of $N$. On the other hand, the fixed field $F=K^{N'}$ provides
an extension $F/k$ with Galois closure $K/k$ ( $K/F$ has order 8, 5 or 3 and $S_5$ has no nontrivial normal subgroup of
order dividing any of these numbers). In the proofs of Propositions 3.2 and 4.11 in \cite{CRV2} , mostly
arguing on solvability of holomorphs of groups of order 15, 24 and 40, respectively, we proved that $F/k$ is not Hopf Galois.
In this way, theorem 2.9 in \cite{KKTU} rules out  all these 44  Hopf Galois types.
We perform a  computation with Magma to check that the
holomorph of $\SL(2,5)$ does not contain $S_5$ as a transitive subgroup and  then we have the following result.

\begin{proposition}
Let $K/k$ be a Galois extension with Galois group $A_5$. Then, the only type of Hopf Galois structures on $K/k$ is $A_5$.
The classical Galois structure realizes this type.

Let $K/k$ be a Galois extension with Galois group $S_5=A_5\rtimes C_2$. Then, the only types of Hopf Galois structures on $K/k$ are $S_5$
and the split one $A_5\times C_2$. The classical Galois structure realizes the first type and the second type is realized as the
induced Hopf Galois structure by an almost classical Hopf Galois structure on $K^{<\tau>}/k$, where $\tau$ denotes
a transposition in $S_5$.
\end{proposition}

\subsection{Galois extensions with Galois group $A_n$ or $S_n$, $n \geq 5$}\label{n}

Let $K/k$ be a Galois extension with Galois group $G=S_n$ or $A_n$, $n \ge 5$. Let us assume that $K/k$ is Hopf Galois of cyclic type. If $N$ is a cyclic group corresponding to this Hopf Galois structure, we can take for $N'$ the unique subgroup of $N$ of index $n$. Then $K^{N'}/k$ should be Hopf Galois (of cyclic type). But we know by \cite{G-P} Corollary 4.8 that a separable extension of degree $n\ge 5$ such that its normal closure has Galois group $A_n$ or $S_n$ is not
Hopf Galois. We have then obtained the following result.

\begin{proposition}\label{f}
Let $K/k$ a Galois extension with Galois group $S_n$ or $A_n$, where $n\ge 5$. Then $K/k$ has no
Hopf Galois structures of cyclic type.
\end{proposition}

Let us note that the results for the alternating group in sections \ref{5} and \ref{n} are special cases of Byott's main result in \cite{B2} where the author proves that a Galois extension $K/k$ with Galois group a non-abelian simple group $G$ has exactly two Hopf Galois structures the Galois one and the classical non-Galois one. The results for the symmetric group complements those in \cite{C-C}, where the authors compute the number of Hopf Galois structures of types $S_n$ and $A_n\times C_2$ on a Galois extension with Galois group $S_n$. Proposition \ref{f} supports a query of Byott in \cite{B3} where he states that we do not have any examples where an extension with nonsolvable Galois group admits a Hopf Galois structure of solvable type.


\begin{thebibliography}{1}
\bibitem{B} N.P. Byott, \emph{Uniqueness of Hopf Galois structure for separable field extensions}. Comm.  Algebra 24 (1996), 3217-3228. Corrigendum, ibid., 3705.
\bibitem{B2} N.P. Byott, \emph{Hopf-Galois structures on field extensions with simple Galois groups}. Bull. London Math. Soc. 36 (2004), 24-29.
\bibitem{B3} N.P. Byott, \emph{Solubility criteria for Hopf-Galois structures}. New York J. Math. 21 (2015), 883-903.
\bibitem{C-C} S. Carnahan, L. N. Childs, \emph{Counting Hopf Galois Structures on Non-Abelian
Galois Field Extensions}. J. Algebra 218 (1999), 81-92.
\bibitem{C-S} S.U. Chase, M. Sweedler, \emph{Hopf Algebras and Galois Theory}. Lecture Notes in Mathematics, Vol. 97, Springer Verlag, 1969.
\bibitem{Ch1} L. N. Childs, \emph{On the Hopf Galois theory for separable field extensions}. Comm. Algebra 17 (1989), 809-825.
\bibitem{Ch2} L. N. Childs, \emph{Taming wild extensions: Hopf algebras and local Galois module theory}, AMS 2000.
\bibitem{CRV} T. Crespo, A. Rio, M. Vela, \emph{On the Galois correspondence theorem in separable Hopf Galois theory}, Publ. Mat. 60 (2016), 221-234.
\bibitem{CRV2} T. Crespo, A. Rio, M. Vela,  \emph{The Hopf Galois property in subfield lattices}, Comm. Algebra 44 (2016), 336-353; Corrigendum, ibid., 3191.
\bibitem{CRV3} T. Crespo, A. Rio, M. Vela, \emph{Induced Hopf Galois structures}, J. Algebra 457 (2016), 312-322.
\bibitem{C-S} T. Crespo, M. Salguero, \emph{An algorithm to determine Hopf Galois structures}, arXiv: 1704.00232.
\bibitem{G-P} C. Greither, B. Pareigis, \emph{Hopf Galois theory for separable field extensions}, J. Algebra 106 (1987), 239-258.
\bibitem{KKTU} A. Koch, T. Kohl, P.J. Truman, R. Underwood, \emph{Normality and short exact sequences of Hopf-Galois structures}, arXiv: 1708.08402.
\bibitem{K} T. Kohl, \emph{Multiple Holomorphs of Dihedral and Quaternionic Groups}, Communications in Algebra, 43 (2015), 4290-4304.
\end{thebibliography}
\end{document}